\documentclass{cmslatex}
\usepackage{latexsym, amssymb, enumerate, amsmath}
\usepackage[dvips]{epsfig}
\usepackage{latexsym,amssymb}

\usepackage{graphics}
\renewcommand{\theequation}{\arabic{section}.\arabic{equation}}

\newcommand{\ens}[1]{\mathbb{#1}}
\newcommand{\ron}[1]{\mathcal{#1}}
\newcommand{\N}{\mathbb{N}}

\newcommand{\R}{\mathbb{R}}
\newcommand{\C}{\mathbb{C}}

\def\cal{\mathcal}

\def\derpar#1#2{\frac{\partial#1}{\partial#2}}

\newtheorem{thm}{Theorem}[section]

\newtheorem{rem}[thm]{Remark}
\newtheorem{rems}[thm]{Remarks}
\begin{document}

\title{Quantitative linearized study of the Boltzmann collision operator 
and applications
\thanks{%{Received date / Revised version date}
          % The correct dates will be entered by the CMS editor}}
}}
          %For each author, make a block with the following four macros:
\author{Cl\'ement Mouhot
\thanks {CEREMADE, Universit\'e Paris IX Dauphine
Place du Mar\'echal de Lattre de Tassigny
75775 PARIS Cedex 16, FRANCE, cmouhot@ceremade.dauphine.fr.}
%\and full name \thanks {address, (email).}
}
          %{Put the URL for your home page here if you have one}

          %Use \thanks statements for acknowledgements of grants and
          %support. They will appear below all the authors' addresses, so be
          %specific about which author is thanking whom:

          \thanks{Support by the European network HYKE, funded by the EC as
contract HPRN-CT-2002-00282, is acknowledged.}

\pagestyle{myheadings} \markboth{Quantitative linearized study of the Boltzmann collision operator}{C. Mouhot}\maketitle

\begin{abstract}
We present recent results \cite{BaMo,Mo:coerc:05,GM:04} 
about the quantitative study of the linearized Boltzmann 
collision operator, and its application 
to the study of the trend to equilibrium for the spatially 
homogeneous Boltzmann equation for hard spheres. 
\end{abstract}

\begin{keywords}
Boltzmann equation, spatially homogeneous,  
linearized Boltzmann collision operator, spectrum, spectral gap,  
explicit, trend to equilibrium, rate of convergence. 
\medskip

{\bf AMS subject classifications.}
76P05 Rarefied gas
flows, Boltzmann equation [See also 82B40, 82C40, 82D05].
\end{keywords}

\section{Introduction}
\setcounter{equation}{0}

In this paper, we present some of our recent works \cite{BaMo,Mo:coerc:05,GM:04} 
(the first one \cite{BaMo} being in collaboration with C\'eline Baranger) on 
the quantitative study of the linearized Boltzmann collision 
operator and its application to the quantitative study of trend to equilibrium. 
This first section shall be devoted to the introduction of the Boltzmann 
equation and the motivation of this study. Section~\ref{sec:clas} 
recalls classical results on the linearized Boltzmann collision operator, 
and Section~\ref{sec:BM} and~\ref{sec:coerc} present some new explicit estimates 
on this operator. Finally Section~\ref{sec:GM} present an application 
of these explicit estimates in the study of explicit rate of convergence 
to equilibrium for the (nonlinear) spatially homogeneous Boltzmann equation. 

\subsection{The Boltzmann equation.}
The {\em Boltzmann equation} describes the behavior of
a dilute gas when the only interactions taken into account are 
binary collisions, by means of an evolution equation 
on the time-dependent particle distribution function in 
the phase space. 
In the case where this distribution function is
assumed to be independent of the position, we obtain 
the {\em spatially homogeneous Boltzmann equation}:
 \begin{equation}\label{eq:base}
 \derpar{f}{t}  = Q(f,f), \quad  v \in \R^N, \quad t \geq 0
 \end{equation}
in dimension $N \ge 2$. In spite of the strong restriction that this assumption 
of spatial homogeneity constitutes, it has proven an interesting 
and inspiring case for studying qualitative properties of the 
Boltzmann equation. In equation~\eqref{eq:base}, $Q$ is the quadratic
Boltzmann collision operator, defined by the bilinear form
 \begin{equation*}\label{eq:collop}
 Q(g,f) = \int _{\R^N \times \ens{S}^{N-1}} B \, (g'_* f' - g_* f) \, dv_* \, d\sigma.
 \end{equation*}
Here we have used the shorthands $f'=f(v')$, $g_*=g(v_*)$ and
$g'_*=g(v'_*)$, where
 \begin{equation*}\label{eq:rel:vit}
 v' = \frac{v+v_*}2 + \frac{|v-v_*|}2 \, \sigma, \qquad
 v'_* = \frac{v+v_*}2 - \frac{|v-v_*|}2 \, \sigma
 \end{equation*}
stand for the pre-collisional velocities of particles which after
collision have velocities $v$ and $v_*$. 
$B$ is the Boltzmann {\em collision kernel} determined by physics 
(related to the cross-section $\Sigma(v-v_*,\sigma)$ 
by the formula $B=|v-v_*| \, \Sigma$). On physical grounds, it is
assumed that $B \geq 0$ and $B=B(|v-v_*|, \cos \theta)$ is a function of $|v-v_*|$ and
$\cos\theta$, where $\theta\in [0,\pi]$ is the {\em deviation angle} 
between $v'-v'_*$ and $v-v_*$, defined by 
$$
\cos \theta = \frac{v-v_*}{|v-v_*|} \cdot \sigma.
$$
\smallskip

Boltzmann's collision operator has the fundamental properties of
conserving mass, momentum and energy
  \begin{equation*}
  \int_{\R^N}Q(f,f) \, \phi(v)\,dv = 0, \quad
  \phi(v)=1,v,|v|^2 %\label{eq:CON}
  \end{equation*}
and satisfying Boltzmann's $H$ theorem, which can be formally written as 
  \begin{equation*} %\label{eq:HTH}
  {\cal D}(f):= - \frac{d}{dt} \int_{\R^N} f \log f \, dv = - \int_{\R^N} Q(f,f)\log(f) \, dv \geq 0.
  \end{equation*}
The $H$ functional $H(f) = \int f \log f$ is the opposite of the entropy of the solution.
Boltzmann's $H$ theorem implies that any equilibrium distribution
function has the form of a {\em Maxwellian distribution}
  \begin{equation*}
  M(\rho,u,T)(v)=\frac{\rho}{(2\pi T)^{N/2}}
  \exp \left( - \frac{\vert u - v \vert^2} {2T} \right), %\label{eq:MAX}
  \end{equation*}
where $\rho,\,u,\,T$ are the density, mean velocity
and temperature of the gas
  \begin{equation*}
  \rho = \int_{\R^N}f(v) \, dv, \quad u =
  \frac{1}{\rho}\int_{\R^N}vf(v) \, dv, \quad T = {1\over{N\rho}}
  \int_{\R^N}\vert u - v \vert^2f(v) \, dv, 
  \end{equation*}
which are determined by the mass, momentum and energy of the initial datum thanks 
to the conservation properties. As a result of the process of entropy production 
pushing towards local equilibrium combined with the constraints of conservation laws, 
solutions are thus expected to converge to a unique Maxwellian equilibrium.  
Up to a normalization we set without restriction $M(v) = e^{-|v|^2}$ 
as the Maxwellian equilibrium, or equivalently $\rho=\pi^{N/2}$, 
$u=0$ and $T=1/2$. 

\subsection{Motivation.}
The relaxation to equilibrium is studied since the works of Boltzmann and 
it is at the core of the kinetic theory. The motivation is to provide 
an analytic basis for the second principle of thermodynamics for a statistical physics 
model of a gas out of equilibrium. 
Indeed Boltzmann's famous $H$ theorem gives an analytic meaning to 
the entropy production process and identifies possible 
equilibrium states. In this context, proving convergence towards  
equilibrium is a fundamental step to justify Boltzmann model, but 
cannot be fully satisfactory as long as it remains based on non-constructive 
arguments. Indeed, as suggested implicitly by Boltzmann when answering 
critics of his theory based on Poincar\'e recurrence Theorem, 
the validity of the Boltzmann equation breaks for very 
large time (see~\cite[Chapter~1, Section~2.5]{Vi:hb} 
for a discussion). It is therefore crucial to obtain quantitative 
informations on the time scale of the convergence, in order to show that 
this time scale is much smaller than the time scale of 
validity of the model. 
Moreover constructive arguments often provide new qualitative insights into the 
model, for instance here they give a better understanding of 
the dependency of the rate of convergence according to the 
collision kernel and the initial datum. 

\subsection{Assumptions on the collision kernel.}
The main physical case of application of this paper is that of 
hard spheres in dimension $N=3$, where (up to a normalization constant) 
$B (|v-v_*|, \cos \theta) = |v-v_*|$. 
More generally we shall make the following decoupling assumption 
on the collision kernel: we assume that $B$ takes the product form
  \begin{equation}\label{eq:prod}
  B(|v-v_*|,\cos \theta) = \Phi (|v-v_*|) \, b(\cos \theta),
  \end{equation}
where $\Phi$ and $b$ are nonnegative functions not identically equal to $0$. 
This decoupling assumption is made for the sake of 
simplicity and could probably be relaxed at the price of 
technical complications. Additional assumptions on the collision kernel shall be given in 
each section. Let us recall nevertheless the fundamental 
class of {\em inverse power-law interaction models}, for which $\Phi(z) = C_\Phi \, |z|^\gamma$ 
for some $\gamma \in (-N,1]$ and $C_\Phi>0$, and $\theta \mapsto b(\cos \theta)$ is not 
integrable at $\theta \sim 0$. The classical vocabulary denotes 
by {\em hard potentials} the case $\gamma >0$, {\em Maxwell molecules} 
the case $\gamma =0$, and {\em soft potentials} the case $\gamma <0$. 
{\em Grad's angular cutoff assumption} (see \cite{Grad58}) means a truncation 
so that $\theta \mapsto b(\cos \theta)$ is integrable on $[0,\pi]$. 
\smallskip

In the case where $B$ is locally integrable, we can 
define the so-called {\em collision frequency} 
  \begin{equation*} \label{eq:colfreq}
  \nu(v) = \int_{\R^N \times \ens{S}^{N-1}} \Phi(|v-v_*|) \, b(\cos \theta) \, M(v_*) \, dv_* \, d\sigma. %= (\Phi * M) (v).
  \end{equation*}
We denote by $\nu_0 \ge 0$ the minimum value of $\nu$. It is easily seen 
that for inverse power-law collision kernels, this minimum value is positive as soon as $\gamma \ge 0$. 

\subsection{Linearization.}
We define the {\em linearized collision operator} $\ron{L}$ by 
  \[ \ron{L}(g) = Q(M+g,M+g) - Q(g,g) = Q(M,g) + Q(g,M). \]
We shall consider this operator on several functional spaces. The classical 
spectral theory of the Boltzmann collision operator, as pioneered by 
Hilbert~\cite{Hilb:EB:12}, is done on $g \in L^2(M^{-1})$. 
However the Cauchy theory for the nonlinear spatially 
homogeneous Boltzmann equation corresponds 
to the functional space $g \in L^1(1+|v|^2)$. We shall therefore introduce 
intermediate functional spaces in order to connect these two theories in the last section.
\smallskip
 
%$L^1(\exp(a |v|^s)$,  
%with $a >0$ and $0<s<2$ to be chosen later, and study the linearized 
%collision operator $\ron{L}$ on these functional spaces.

Let us recall some spectral theory's notions. 
Let us consider a linear unbounded operator $T:\cal{B} \to \cal{B}$ on the 
Banach space $\cal{B}$, defined on a dense domain $\mbox{Dom}(T) \subset \cal{B}$. 
Then we adopt the following notations and definitions:
  \begin{itemize}
  \item we denote by $N(T) \subset \cal{B}$ the {\em null space} of $T$; 
  %\item we denote by $R(T) \subset \cal{B}$ the {\em range} of $T$; 
  \item $T$ is said to be {\em closed} if its graph is closed 
  in ${\cal B} \times {\cal B}$.
  \end{itemize}
\smallskip

In the following definitions, $T$ is assumed to be closed:
%, and it is said to be {\em closable} if any sequence $x_n$ 
%of ${\cal B} \cap D(T)$ going to $0$ in ${\cal B}$ and such that $Tx_n$ 
%converges in ${\cal B}$ satisfies $Tx_n \to 0$.  
  \begin{itemize}
  \item the {\em resolvent set} of $T$ denotes 
        the set of complex numbers $\xi$ such that $T-\xi$ is 
        bijective from $\mbox{Dom}(T)$ to $\cal{B}$ and the 
        inverse linear operator $(T-\xi)^{-1}$, 
        defined on $\cal{B}$, is bounded (see \cite[Chapter~3, Section 5]{Kato});
  \item we denote by $\Sigma(T) \subset \C$ the {\em spectrum} of $T$, that is 
        the complementary set of the resolvent set of $T$ in $\C$;
  \item an {\em eigenvalue} is a complex number $\xi \in \C$ such that 
        $N(T-\xi)$ is not reduced to $\{0\}$;
  \item we denote $\Sigma_d(T) \subset \Sigma(T)$ the {\em discrete spectrum} of $T$, 
        {\em i.e.}, the set of {\em discrete eigenvalues}, that is the eigenvalues 
        isolated in the spectrum and with finite multiplicity ({\em i.e.}, such that the spectral projection  
        associated with this eigenvalue has finite dimension, see \cite[Chapter~3, Section 6]{Kato});
  \item for a given discrete eigenvalue $\xi$, we shall call the {\em eigenspace} 
        of $\xi$ the range of the spectral projection associated with $\xi$;  
  \item we denote $\Sigma_e(T) \subset \Sigma(T)$ the {\em essential spectrum} of $T$ 
        defined by $\Sigma_e(T) = \Sigma(T) \setminus \Sigma_d(T)$;
  \item when $\Sigma(T) \subset \R_-$, we say that $T$ has a {\em spectral gap} when the distance 
  between $0$ and $\Sigma(T) \setminus \{0\}$ is positive, and the spectral gap denotes this distance;
  \item the operator $T$ is said to be {\em sectorial} when 
  $\Sigma(T) \subset \{|\arg(\xi - \lambda)| \ge \pi/2 + w \}$ for some $w > 0$ and some $\lambda \in \R$, 
  and its resolvent satisfies the control: for any $\varepsilon >0$, there is $M_\varepsilon >0$ 
  such that $\|(T-\xi)^{-1}\| \le M_\varepsilon /|\xi - \lambda|$ for 
  $\xi \in \{|\arg(\xi - \lambda)| \le \pi/2 + w - \varepsilon \}$ (see \cite[Chapter~9, Section~1]{Kato}). 
  \end{itemize}
\smallskip

\subsection{Notation.} 
In the sequel we shall denote 
$\langle \cdot \rangle = \sqrt{1 + |\cdot|^2}$. 
For any Borel function $w:\R^N \to \R_+$, we define the 
weighted Lebesgue space $L^p (w)$ on $\R^N$ ($p \in [1,+\infty]$), by the norm 
 \[ \| f \|_{L^p(w)} = \left[ \int_{\R^N} |f (v)|^p \, w(v) \, dv \right]^{1/p} \] 
if $p < +\infty$ and 
 \[ \| f \|_{L^\infty (w)} = \sup_{v \in \R^N} |f (v)| \, w(v) \]
when $p = +\infty$. 
The weighted Sobolev spaces $W^{k,p} (w)$ 
($p \in [1,+\infty]$ and $k \in \N$) are defined by the norm 
 \[ \| f \|_{W^{k,p} (w)} = 
       \left[ \sum_{|s| \le k} \|\partial^s f \|^p _{L^p(w)} \right]^{1/p} \]
with the notation $H^k (w) = W^{k,2} (w)$. 
%In the sequel we shall denote by $\| \cdot \|$ indifferently 
%the norm of an element of a Banach space or the usual 
%operator norm on this Banach space, and we shall denote 
%by $C$ various positive constants independent of the collision kernel. 

\section{The classical linearized theory, from Hilbert to Grad}\label{sec:clas}
\setcounter{equation}{0}

%Let us denote $\Hil=L^2(M^{-1})$ and $\Hil(\nu^2)=L^2(\nu^2 M^{-1})$.
It is well-known from the classical theory of the linearized operator  
(see \cite{Grad63} or \cite[Chapter 7, Section 1]{CIP}) that 
the operator $\ron{L}$ on $L^2(M^{-1})$ with domain $L^2(\nu^2 M^{-1})$ is closed and 
self-adjoint, and satisfies
 \begin{multline*}
 \langle g,\ron{L} g \rangle_{L^2(M^{-1})} = \int_{\R^N} g (\ron{L}g) \, M^{-1} \, dv = 
 - \frac{1}{4}\int_{\R^N \times 
 \R^N \times \ens{S}^{N-1}} \Phi(|v-v_*|)\, b(\cos \theta) \\ 
 \Bigg[ \left(\frac{g}{M} \right)^{'}_{*}+ \left(\frac{g}{M} \right)^{'} 
        - \left(\frac{g}{M} \right)_{*} - \left(\frac{g}{M} \right) \Bigg]^2 
 M\, M_{*} \, dv \, dv_{*} \, d\sigma \le 0.
 \end{multline*}
%, together with the self-adjointness, NON cette inegalite suffit a l'impliquer
%  en fait ca implique meme l'auto-adjonction vu que le numerical range 
%  est alors inclus dans R...
This implies that its spectrum is included in $\R_-$. Its null space is 
  \begin{equation*}\label{noyauL}
  N_{L^2(M^{-1})}(\ron{L}) = \mbox{Span} \left\{ M, v_1 M,\dots,v_N M, |v|^2 M \right\}.
  \end{equation*}
These two properties correspond to the linearization of 
Boltzmann's $H$ theorem. 
\smallskip

Let us denote by $D (g)=-\langle g,\ron{L} g \rangle_{L^2(M^{-1})}$ the Dirichlet form for 
$-\ron{L}$. The existence of a spectral gap $\lambda>0$ can be 
written as 
 \begin{equation*} \label{eq:sg}
 \forall \, g \in L^2(M^{-1}), \ g \bot N_{L^2(M^{-1})}(\ron{L}), \hspace{0.3cm}
          D(g) \ge \lambda \, \|g\|^2 _{L^2(M^{-1})}. 
 \end{equation*}
 \smallskip

The first step in the study of the linearized collision operator $\ron{L}$ was done by Hilbert (on the 
hard spheres model), who introduced a decomposition $\ron{L} = \ron{K} - \nu$ 
of $\ron{L}$ between a non-local part $\ron{K}$ (compact in $L^2(M^{-1})$) and 
a multiplication part by the collision frequency $\nu$. This result 
was used by Hilbert for the construction of the so-called 
Hilbert expansion {\em via} the Fredholm alternative. This decomposition 
was then used by Carleman~\cite{Carl57} to show the existence of a spectral gap for 
the hard spheres model. 
Then Grad~\cite{Grad63} generalized this decomposition to hard potentials with angular cutoff or 
Maxwell molecules with angular cutoff, and obtained the existence of a spectral gap 
for these interactions. 
The proof relied on Weyl's Theorem for self-adjoint operators, 
which asserts that the essential spectrum is stable under (relatively) 
compact perturbation. Hence the essential spectrum of $\ron{L}$ is given 
by the one of $-\nu$ (which lies a positive distance away from $0$ 
when $\gamma \ge 0$). Since the operator is non-positive, 
it has only discrete eigenvalues in $(-\nu_0,0]$, possibly accumulating only at $-\nu_0$. 
\smallskip

It was also already observed that the Dirichlet form is monotonous according 
to $B$, which enables to prove the existence of a spectral gap 
for any collision kernel $B$ controlled from below by some 
cutoff hard potentials collision kernel. Moreover 
as discussed in~\cite[Chapter~4, Section~6]{Ce88}, it was proved in~\cite{KuWi67} 
that $\ron{L}$ has an infinite number of discrete negative eigenvalues in the interval 
$(-\nu_0,0)$, which implies that the spectral gap $\lambda$ satisfies $0 < \lambda < \nu_0$. 
In fact the proof in~\cite{KuWi67} was done for hard spheres, but the 
argument applies to any cutoff hard potential collision kernel as well 
(see~\cite[Chapter~4, Section~6]{Ce88}). 
\smallskip

Before going into the next section, let us note that this classical 
theory of the linearized collision operator made an important 
breakthrough: it showed exponential convergence in some linearized 
settings, and motivated subsequent development of the perturbative theory of 
the Boltzmann equation, see \cite{Ukai74} for instance where 
the first smooth global solutions for the spatially inhomogeneous Boltzmann equation 
were constructed. However, its main drawback was that it was 
unable to provide any information on the size of the spectral gap. 

\section{Explicit spectral gap estimates}\label{sec:BM}
\setcounter{equation}{0}

At about the same time as the work of Grad, another theory appears in the study of the 
Boltzmann equation. In the particular case of Maxwell molecules ($\gamma=0$), 
it has been obtained (see \cite{WCUh52,WCUh70}) 
a complete and explicit study of the spectrum of $\ron{L}$ in $L^2(M^{-1})$, 
by symmetry arguments. Indeed in the particular case where the collision 
kernel is independent on the modulus of the relative velocity, the 
collision operator lends itself to an extensive Fourier transform 
analysis, as first noticed by Bobylev, see \cite{Boby88}. 
\smallskip

However, whenever $\Phi$ depends on the modulus of the relative velocity and 
the collision frequency is bounded from below by some positive number 
(for instance for the important physical case of hard spheres), 
the spectral gap was shown to exist by Grad's argument, but no quantitative information was provided, 
in particular on its size. 
%For the important physical case of hard spheres, 
%no information was available on this spectral gap. 
In the work \cite{BaMo}, we propose a new method to show the existence 
of a spectral gap for any hard potentials (in the generalized sense~\eqref{eq:hypPhi}), 
which is geometrical and based on a physical argument. It gives explicit
estimates and deals with the whole operator, with or without angular cutoff.  
\smallskip

Let us write down the assumptions for the collision kernel $B$ (in addition 
to the decoupling assumption~\eqref{eq:prod}):
 \begin{itemize}
 \item The kinetic part $\Phi$ is bounded from below at infinity, {\em i.e.},  
  \begin{equation}\label{eq:hypPhi}
  \exists \, R \ge 0, \ c_\Phi >0 \ \  | \ \ \forall \, r \ge R, \ \Phi(r) \ge c_\Phi.
  \end{equation}
 This assumption holds for hard potentials (and hard spheres).
 \item The angular part $b$ satisfies
 \begin{equation}\label{eq:hypb}
  c_b = \inf_{\sigma_1, \sigma_2 \in \ens{S}^{N-1}} 
  \int_{\sigma_3 \in \ens{S}^{N-1}} \min \{  b(\sigma_1 \cdot \sigma_3), 
  b(\sigma_2 \cdot \sigma_3) \} \, d\sigma_3 > 0. 
  \end{equation}
 This assumption covers all the physical cases.
 \end{itemize}
\medskip

The main theorem of this work is 
\smallskip

 \begin{thm} \label{theo:Bolt}
 Under the assumptions~\eqref{eq:prod}, \eqref{eq:hypPhi}, \eqref{eq:hypb}, 
 the Dirichlet form $D$ (in $L^2(M^{-1})$) of the linearized 
 collision operator with collision kernel $B = \Phi \, b$   
 satisfies
  \begin{equation*}\label{eq:reducB}
   \forall \, g \in L^2(M^{-1}), \ g \bot N_{L^2(M^{-1})}(\ron{L}), \hspace{0.3cm} D (g) \ge C _{\Phi,b} \, D _0 (g),
  \end{equation*}
 where $D _0 (g)$ stands for the Dirichlet form (in $L^2(M^{-1})$) of the 
 linearized collision operator with collision kernel $B_0 \equiv 1$ and
  \begin{equation*}
  C _{\Phi,b} = \left( \frac{c_\Phi \, c_b \, e^{-4 R^2}}{32 \left|\ens{S}^{N-1}\right|} \right)  
  \end{equation*}
 with $R$, $c_\Phi$, $c_b$ being defined in~\eqref{eq:hypPhi}, \eqref{eq:hypb}.
 
 As a consequence we deduce that
  \begin{equation*}\label{eq:spB}
  \forall \, g \in L^2(M^{-1}), \ g \bot N_{L^2(M^{-1})}(\ron{L}), \hspace{0.3cm} 
  D (g) \ge C _{\Phi, b} \, |\lambda _0| \, \| g \|^2 _{L^2(M^{-1})}.
  \end{equation*}
 Here $\lambda _0$ is the spectral gap (that is the modulus of 
 the first non-zero eigenvalue) in $L^2(M^{-1})$ of the 
 linearized Boltzmann operator with $B_0 \equiv 1$, which equals in dimension $3$ 
 (for density $\pi^{3/2}$, momentum $0$ and temperature $1/2$, see~\cite{Boby88})
  \begin{equation*} 
  \lambda _0 = \pi \, \left( \frac{\pi}2 \right)^{3/2} \, \int_0 ^\pi \sin ^3 \theta \, d\theta 
  = \left( \frac{\pi}2 \right)^{3/2} \, \frac{4 \pi}{3}. 
  \end{equation*}
 \end{thm}
\medskip

\begin{rem} Let us mention that in the same work, similar results are derived for the 
linearized Landau operator for hard potentials by grazing collision limit. 
\end{rem}
\medskip

As an application of Theorem~\ref{theo:Bolt}, let us give explicit  
controls for the spectral gap $\lambda = \lambda_{b,\gamma}(\rho,u,T)$ of the 
Boltzmann linearized operator with $b \ge 1$ 
and $\Phi(z) = |z|^\gamma$, $\gamma >0$, in dimension $3$, and 
for a steady state with density $\rho \ge 0$, momentum 
$u \in \R^3$, and temperature $T >0$. By scaling and translation 
arguments, it is straightforward that 
$$
\lambda_{b,\gamma}(\rho,u,T) = \rho \, T^{\left(\frac{N+\gamma}2\right)} \, \lambda_{b,\gamma}(1,0,1). 
$$
\smallskip
 
Then $c_b \ge |\ens{S}^2|$ 
and for any given $R$ we can take $c_\Phi = R^\gamma$. Thus we get 
 \[ \lambda_{b,\gamma}(1,0,1) \ge 2^{\gamma/2} \, \left( \frac{R^\gamma \, e^{-4 R^2}}{32} \right) \, \frac{4 \pi}{3}  \] 
for any $R >0$ (the factor $2^{\gamma/2}$ comes from the fact that $\lambda_0$ was computed for 
density $\pi^{3/2}$, momentum $0$ and temperature $1/2$). An easy computation leads to the lower bound 
 \[ \lambda_{b,\gamma}(1,0,1) \ge 2^{\gamma/2} \, \frac{\pi \, (\gamma/8)^{\gamma/2} \, e^{-\gamma/2}}{24} \]
by optimizing the free parameter $R$. For instance, 
in the hard spheres case ($b=1$ and $\gamma=1$), we obtain a lower bound 
$\lambda_{1,1}(1,0,1) \ge c \approx 0.04$. 
\medskip

The idea of the proof is to work directly on the Dirichlet form 
(without requiring to some decomposition of the collision operator and 
perturbative theory anymore) and to reduce the case of hard potentials 
(in the generalized sense~\eqref{eq:hypPhi}) 
to the Maxwellian case. The difficulty is to deal with the 
cancellations of the collision kernel, and 
most importantly those of the kinetic part $\Phi$ on the diagonal $v = v_*$. 
\smallskip

A first ingredient is the following inequality (which is a corollary
of~\cite[Theorem 2.4]{CMCV:gran:01})
 \begin{multline}\label{eq:CMCV}
 \int_{\R^N} \int_{\R^N} | \varphi (x) - \varphi(y)|^2 \, |x-y|^\gamma \, M(x) \, M(y) \, dx \, dy \\
 \ge K_\gamma \,  \int_{\R^N} \int_{\R^N} | \varphi (x) - \varphi(y)|^2  \, M(x) \, M(y) \, dx \, dy
 \end{multline}
for $\gamma \ge 0$, $\varphi$ some function, and 
 \begin{equation*}
 K_\gamma = \frac{1}{4 \int_{\R^N} M} \ \inf_{x,y \in \R^N} \int_{\R^N} 
 \min \left\{ |x-z|^\gamma , |z-y|^\gamma \right\} \, M(z) \, dz .
 \end{equation*}
\smallskip

\begin{figure}
 \centering
 \includegraphics[width=10cm]{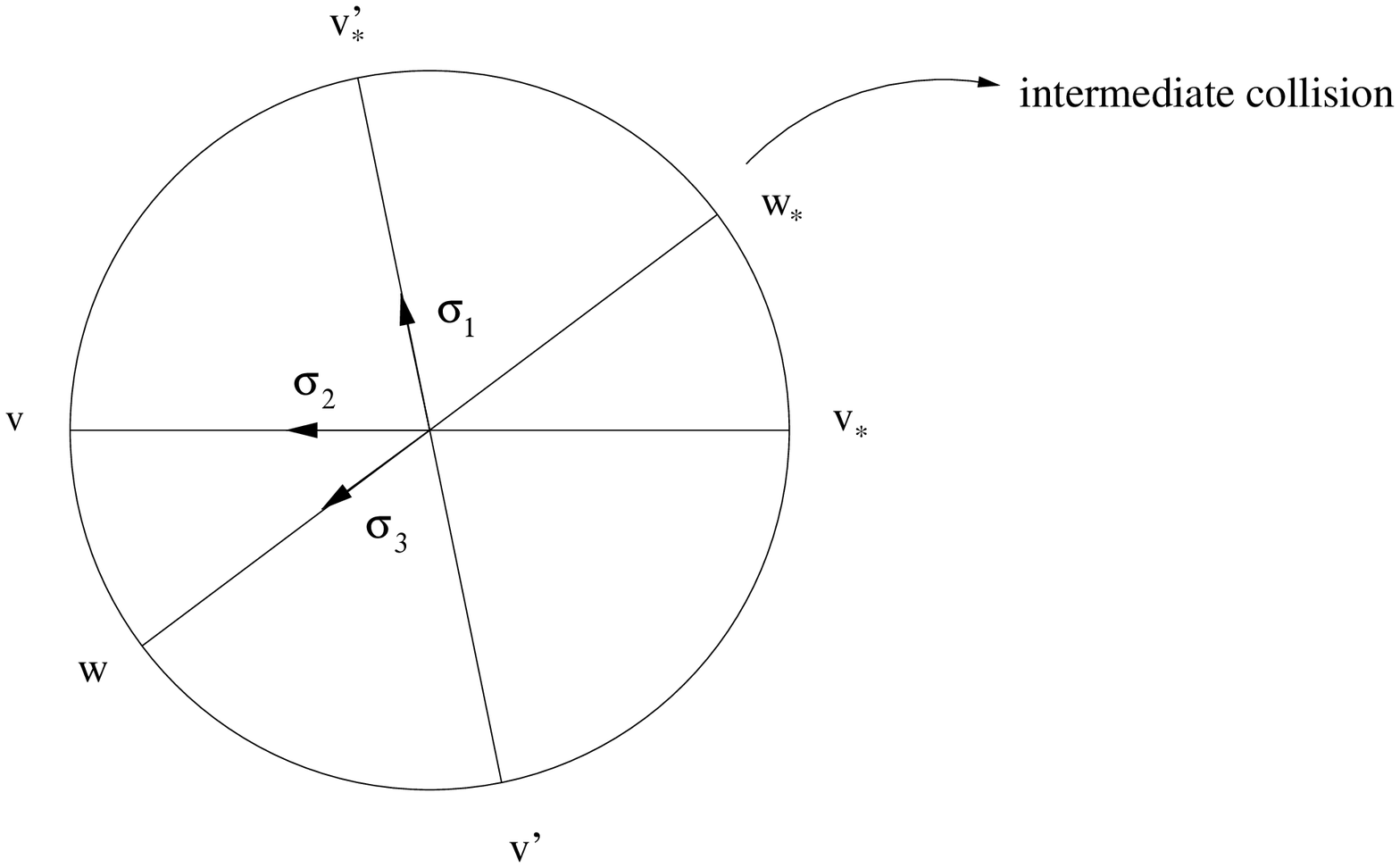}
 \caption{Introduction of an intermediate collision to treat cancellations of $b$}\label{fig:ang}
 \end{figure}

\begin{figure}
 \centering
 \includegraphics[width=10cm]{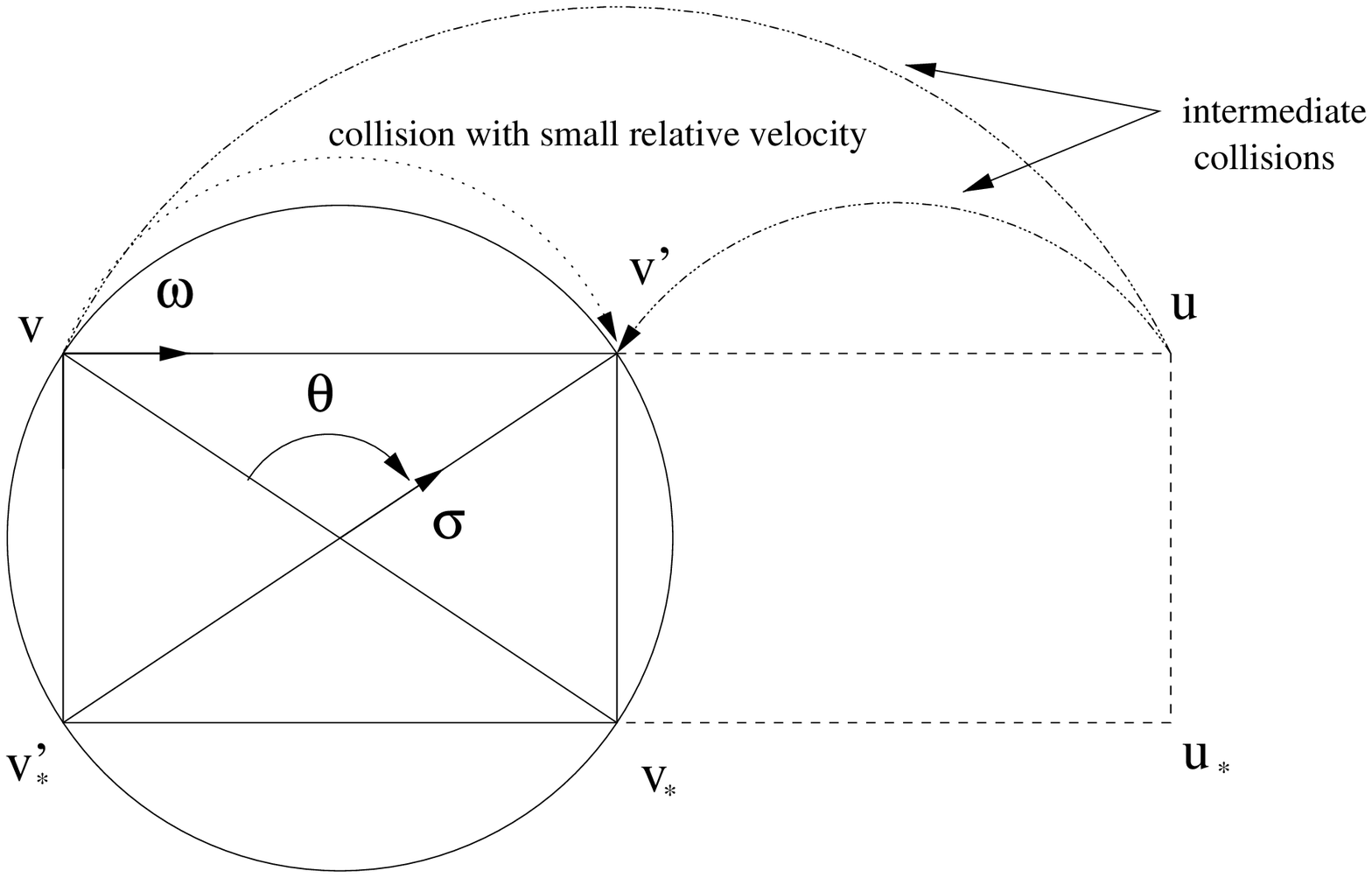}
 \caption{Introduction of an intermediate collision to treat cancellations of $\Phi(z)$ at $z \sim 0$}\label{fig:2D}
 \end{figure}

The proof of~\eqref{eq:CMCV} relies strongly on the existence 
of a ``triangular inequality'' for the function $F(x,y)=|\varphi(x) - \varphi(y)|^2$ integrated, 
which follows immediately here from 
 \[ F(x,y) \le 2 F(x,z) + 2 F(z,y). \]
All the point is thus to adopt suitable representations of the linearized 
collision operator, which contain some integrals of this form, and to estimate these integrals 
using some ``triangular inequality'' 
for the integrated function. Cancellations for $b$ are dealt with by considering 
the integral according to the variables $x=(v-v_*)/|v-v_*|$ and $y=(v'-v'_*)/|v'-v'_*|$ 
on the sphere (see Figure~\ref{fig:ang} where the intermediate variable 
is $z=(w-w_*)/|w-w_*|$).  
Cancellations of $\Phi$ are dealt with by considering the integration 
according to the variables $x=v \cdot (v'-v)/|v'-v|$ and $y = v' \cdot (v'-v)/|v'-v|$ 
on $\R$ (see Figure~\ref{fig:2D} where the intermediate variable is 
$z=u \cdot (v'-v)/|v'-v|$). 
\smallskip

Physically this corresponds to replacing in the Dirichlet form 
the collisions which cancel the collision kernel by some sequences of collisions for which 
the collision kernel is controlled from below. More precisely for instance to treat collisions 
which cancel $\Phi$ (small relative velocity), we replace each collision which has a small relative 
velocity by a sequence (two indeed) of intermediate collisions 
with large relative velocities in order to use \eqref{eq:hypPhi} (see Figure~\ref{fig:2D}). 

\section{Explicit coercivity estimates}\label{sec:coerc}
\setcounter{equation}{0}

In this section, we shall present the results of the work \cite{Mo:coerc:05}, which 
generalize the explicit spectral gap estimates into explicit coercivity estimates 
for any inverse power-law interactions. 
\smallskip

We assume that the collision kernel $B$ satisfies~\eqref{eq:prod} and:
 \begin{itemize}
 \item The kinetic part $\Phi$ is bounded from below by a power-law:
   \begin{equation}\label{eq:hypPhibis}
   \forall \, r \ge 0, \ \ \ \Phi(r) \ge C_\Phi r^\gamma
   \end{equation}
 where $\gamma \in (-N,1]$ and $C_\Phi>0$ is some constant. 
 Collision kernels deriving from interaction potentials behaving like inverse-power laws 
 satisfy this assumption, as well as hard spheres collision kernels. 
 \item The angular part $b$ satisfies~\eqref{eq:hypb}. 
 Moreover, in order to obtain regularity estimates 
 when the collision kernel is not locally integrable, we shall assume the more accurate 
 control from below
   \begin{equation}\label{eq:hypbnc}
   \forall \, \theta \in (0,\pi], \ \ \ b(\cos \theta) \ge \frac{c_b}{\theta^{N-1+\alpha}} 
   \end{equation}
 for some constant $c_b >0$ and $\alpha \in [0,2)$ (note that  
 assumption~\eqref{eq:hypbnc} implies straightforwardly assumption~\eqref{eq:hypb}). 
 The goal of this control is to measure the strength of the angular singularity, which is related 
 to the regularity properties of the collision operator (as already noticed in \cite{ADVW} for instance). 
 \end{itemize}
\smallskip

It was proved in~\cite{Cafl80} that the linearized collision operator 
for soft potentials ($\gamma <0$) with cutoff has no spectral gap. 
But if one allows a loss on the algebraic weight of the norm, 
it was proved in~\cite{GoPo86} a ``degenerated spectral gap'' result 
of the form: 
  \begin{equation}\label{GP}
  \forall \, g \in L^2(M^{-1}), \ g \bot N_{L^2(M^{-1})}(\ron{L}), \hspace{0.3cm} 
  D (g) \ge C \, \big\| g \, \langle v \rangle^{\gamma/2} \big\|^2 _{L^2(M^{-1})}  
  \end{equation} 
where $\gamma <0$ is the exponent in~\eqref{eq:hypPhibis}. The proof was based on 
inequalities proved in~\cite{Cafl80} together with Weyl's Theorem and it 
gave no informations on $C$. 
\smallskip

In the work \cite{Mo:coerc:05}, we extend and complete the works~\cite{GoPo86} and \cite{BaMo} by 
  \begin{itemize}
  \item giving a constructive proof of estimate~\eqref{GP} for soft potentials (with or without cutoff);
  \item extending it to hard potentials ($\gamma>0$)  
  (note that for hard potentials this estimate is stronger than the usual spectral gap estimate); 
  \item giving a coercivity result in local Sobolev spaces for the linearized Boltzmann operator  
  with a non locally integrable collision kernel, and discussing the consequence on its spectrum.
  \end{itemize}
\smallskip

In the following theorem, 
$H^{\alpha/2} _{\mbox{{\scriptsize loc}}}$ denotes the space of functions whose restriction to any 
compact set $K$ of $\R^N$ belongs to 
$H^{\alpha/2}(K) = \{ h \in L^2(K) \mbox{ s. t. } (1-\Delta)^{-\alpha/4}h \in L^2(K) \}$ 
(here $L^2(K)$ denotes space of functions square integrable on $K$).  
\medskip

The main result of this work is:
\smallskip

 \begin{thm}\label{theo:Bolt:coerc}
 Under the assumptions~\eqref{eq:prod}, \eqref{eq:hypPhibis}, \eqref{eq:hypb}, 
 the linearized Boltzmann operator $\ron{L}$ with collision kernel $B = \Phi \, b$ 
 satisfies
   \begin{equation}\label{coBo}
    \forall \, g \in L^2(M^{-1}), \ g \bot N_{L^2(M^{-1})}(\ron{L}), \hspace{0.3cm} 
   D (g) \ge C _{\gamma} \, \big\| g \, \langle v \rangle ^{\gamma/2} \big\|_{L^2(M^{-1})} ^2
   \end{equation}
 where $C_\gamma$ is an explicit constant depending only on 
 $\gamma$, $C_\Phi$, $c_b$, and the dimension $N$. 

 When moreover the collision kernel is not locally integrable and $b$ satisfies \eqref{eq:hypbnc}, 
 $\ron{L}$ satisfies \eqref{coBo} and 
   \begin{equation}\label{coBonc}
    \forall \, g \in L^2(M^{-1}), \ g \bot N_{L^2(M^{-1})}(\ron{L}), \hspace{0.3cm} 
    D (g) \ge C _{\gamma,\alpha} \, \|g \|_{H^{\alpha/2} _{\mbox{{\scriptsize {\em loc}}}}} ^2
   \end{equation}
 where $C_{\gamma,\alpha}$ is an explicit constant depending only on 
 $\gamma$, $\alpha$, $C_\Phi$, $C_b$, $c_b$ and the dimension $N$. 
 \end{thm}
\medskip

\begin{rems} 

1. When the collision kernel is locally integrable, the collision 
frequency $\nu$ is finite, and the estimate~\eqref{coBo} can be written in the following form: 
  \[  \forall \, g \in L^2(M^{-1}), \ g \bot N_{L^2(M^{-1})}(\ron{L}), \hspace{0.3cm}  
      D (g) \ge \bar C \, \| g \, \nu \|_{L^2(M^{-1})} ^2 \]
for some explicit constant $\bar C >0$. 
%where $\nu$ is the collision frequency. 
\smallskip

2. When the collision kernel is not locally integrable and $b$ satisfies~\eqref{eq:hypbnc}, 
a natural conjecture would be that the estimate~\eqref{coBonc} improves into: 
  \[  \forall \, g \in L^2(M^{-1}), \ g \bot N_{L^2(M^{-1})}(\ron{L}), \hspace{0.3cm} 
      D (g) \ge C _{\gamma,\alpha} \, \big\| g \, \langle v \rangle ^{\gamma/2} \big\|_{H^{\alpha/2}(M^{-1})} ^2. \]
%where $\alpha \in [0,2)$ is the order of angular singularity, defined in \eqref{eq:hypbnc}, 
%and $H^{\alpha/2} (M^{-1})$ is the Sobolev space defined by 
%$H^{\alpha/2} (M^{-1}) = \{ g \in L^2(M^{-1}) \mbox{ s. t. } (1-\Delta)^{-\alpha/4}g \in L^2(M^{-1}) \}$.  
We were not able to obtain this coercivity estimate, however we give its consequence in terms of local regularity. 
\smallskip

3. When $\gamma>0$ and $\alpha>0$, one can easily 
deduce from Theorem~\ref{theo:Bolt:coerc} that the operator $\ron{L}$ has 
compact resolvent (see~\cite{Mo:coerc:05}), which implies that its spectrum is purely discrete in this case 
(note that it is also true when $\gamma=0$ and $\alpha \ge 0$ by explicit diagonalization). 
\smallskip

4. Similar results have been derived in the same paper \cite{Mo:coerc:05} 
for the linearized Landau collision operator. We also mention the forthcoming work~\cite{MS} 
where some new spectral gap estimates improving on~\eqref{coBo} in the non-integrable case are obtained. 
\end{rems}
\medskip

In the case of hard potentials, the idea of the proof is to decompose the operator 
between a part satisfying the desired 
coercivity estimate and a bounded part, and use the spectral gap estimates. 
This argument is reminiscent of an argument of Grad \cite[Section~5]{Grad63} used 
to study the decrease of the eivenvectors of the linearized 
Boltzmann operator for hard potentials, and it was already noticed 
in~\cite{BCN:86} for instance. Nevertheless it is the first 
time that it is used to obtain explicit estimates (thanks to 
the results in~\cite{BaMo} presented in the previous section). 
\smallskip

For soft potentials we use a dyadic decomposition the Dirichlet form 
according to the modulus of the relative velocity. Combined with 
technical estimates on the non-local part of the linearized collision operators 
and the spectral gap estimates from the Maxwell case, 
it enables to reconstruct a lower bound with the appropriate weight. 
\smallskip

Finally the proof of the coercivity estimates in local Sobolev spaces for the linearized 
Boltzmann operator with a non locally integrable collision kernel 
is inspired by the previous works~\cite{Lions98,Vi:99,ADVW} on the 
regularizing properties of the nonlinear collision operator  
and by the structure of the linearized 
Landau operator in the grazing collision limit. The main task is to 
identity a ``diffusive part'', which shares coercivity properties 
similar to the one of a fractional Laplacian, and a remaining bounded part. 

%on the full nonlinear collision operator,  
%and by structure of the linearized Landau operator. Indeed the 
%suitable decomposition of $L^{{\ron B}}$ for non locally integrable 
%collision kernels (for which the usual Grad's splitting does not 
%make sense anymore) is directly readable on the 
%linearized Landau operator: the part which becomes the diffusion part 
%in the grazing collision limit is the part which enjoys a coercivity 
%property in Sobolev spaces, and the part which becomes the bounded 
%part in the grazing collision limit is the part which is bounded thanks 
%to the cancellation lemmas. 

\section{Applications to the quantitative study of rate of convergence: 
connection to the nonlinear theory}\label{sec:GM}

We present in this section the work \cite{GM:04}, which 
is devoted to the study of the asymptotic behavior 
of solutions to the spatially homogeneous Boltzmann equation for hard potentials with cutoff. 
On one hand it was proved by Arkeryd~\cite{Ar88} by non-constructive arguments that 
spatially homogeneous solutions (with finite mass and energy) of the Boltzmann 
equation for hard spheres converge towards equilibrium with exponential rate, with no 
information on the rate of convergence and the constants (in fact the proof in this paper 
required some moment assumptions, but the latter can be relaxed with the results 
about appearance and propagation of moments, as can be found in \cite{Wenn:momt:97}). 
On the other hand it was proved in~\cite{MoVi} a quantitative convergence result 
with rate $O(t^{-\infty})$ for these solutions. 
\smallskip

In the work \cite{GM:04} we 
improve and fill the gap between these results by 
 \begin{itemize}
 \item showing exponential convergence towards equilibrium by constructive arguments (with explicit rate 
 and constants);
 \item showing that the spectrum of the linearized 
 collision operator in the narrow space $L^2(M^{-1})$ %($M$ is the equilibrium) 
 dictates the asymptotic behavior of the solution 
 in a much more general setting, as was conjectured in~\cite{CaGaTo} on the basis of the study of the Maxwell case. 
 \end{itemize}
\smallskip

We assume that we deal with {\em hard spheres} collision kernel:
 \begin{equation}\label{eq:HS}
 B(|v-v_*|, \cos \theta) = C_B \, |v-v_*|
 \end{equation}
for some $C_B >0$. The results of \cite{GM:04} apply more generally to hard potentials 
with cutoff collision kernels (under some suitable technical assumptions). 
Let us recall that under this assumption on $B$, it has been shown in 
\cite{Ar72,MW} that for any initial datum $0 \le f_0 \in L^1(1+|v|^2)$, there is 
a unique global nonnegative solution in $L^1(1+|v|^2)$ preserving 
mass, momentum and energy. 
\smallskip

Let us review briefly existing results in the study of convergence to equilibrium 
for the spatially homogeneous Boltzmann equation for hard spheres. 
On the basis of the $H$ theorem and suitable {\em a priori}  
estimates, various authors gave results of 
$L^1$ convergence to equilibrium by compactness arguments 
for the spatially homogeneous Boltzmann equation 
with hard potentials and angular cutoff  
(for instance Carleman~\cite{Carl32}, Arkeryd~\cite{Ar72}, etc.). 
These results provide no information at all on the rate of convergence. 
Then, on the basis of these compactness results and the classical linearized 
theory presented in Section~\ref{sec:clas}, 
Arkeryd gave in~\cite{Ar88} the first (non-constructive) proof of 
exponential convergence in $L^1$ for the spatially homogeneous 
Boltzmann equation with hard potentials and angular cutoff. 
His result was generalized to $L^p$ spaces ($1\le p < +\infty$) by Wennberg~\cite{We93}. 
\smallskip

At the beginning of the nineties, several difficulties had still to be 
overcome in order to get a quantitative result of exponential convergence:
 \begin{itemize}
 \item[(i)] No estimate was available for the spectral gap in $L^2(M^{-1})$ 
 for hard spheres. 
 \item[(ii)] There is no known {\em a priori} estimate for the nonlinear 
 problem showing propagation (let alone appearance) of the norm $L^2(M^{-1})$. 
 Matching results obtained in this space and 
 the physical space $L^1(1+|v|^2)$ is one the main difficulties, and 
 was treated in~\cite{Ar88} by a non-constructive argument. 
 \item[(iii)] Finally, any estimate deduced from a linearization 
 argument is valid only in a neighborhood of the equilibrium, and 
 the use of compactness arguments to deduce that the solution 
 enters this neighborhood (as e.g. in~\cite{Ar88}) would 
 prevent any hope of obtaining explicit estimate. 
 \end{itemize} 
\smallskip

\begin{rem}
Note that in the Maxwell case with cutoff, 
all these difficulties have been solved. We have already mentioned 
that in this case, Wang-Chang and Uhlenbeck~\cite{WCUh52,WCUh70} and then Bobylev~\cite{Boby88} 
were able to obtain a complete and explicit diagonalization 
of the linearized collision operator in $L^2(M^{-1})$. Then specific 
metric well suited to the collision operator for Maxwell molecules allowed to 
achieve the goals sketched at the beginning of this section 
(under additional assumptions on the initial datum), see~\cite{CaGaTo} 
and~\cite{CaLu}. However it seems that the proofs are strongly restricted 
to the Maxwellian case. 
\end{rem}
\medskip

In order to solve the point~(iii), quantitative estimates in the large 
have been developed since the beginning of the nineties, 
directly on the nonlinear equation, by relating the entropy 
production functional to the relative entropy:
see~\cite{CC92,CC94,ToVi99,Vi03,MoVi}. 
The latter paper states, for hard spheres, 
quantitative convergence towards equilibrium 
with rate $O(t^{-\infty})$ for solutions in $L^1(1+|v|^2)$. 
Unfortunately it was proved in~\cite{BobyCerc:conjCer:99} that 
Cercignani's conjecture on entropy production is false  
in the case of hard potentials, even for smooth solutions with strong algebraic decay 
(although ``almost true'', in the sense of~\cite{Vi03}). It means that 
one cannot establish in this setting a {\em linear} inequality relating 
the entropy production functional and the relative entropy, 
which would yield exponential convergence directly on the nonlinear equation.
\smallskip

Point~(i) was solved in~\cite{BaMo} (see Section~\ref{sec:BM}). 
\smallskip

In order to solve the remaining obstacle of point~(ii), 
the strategy of \cite{GM:04} is to prove explicit linearized estimates  
of convergence to equilibrium in the space $L^1(\exp ( a|v|^s))$ with 
$a>0$ and $0<s<\gamma/2$, on which we have explicit results of appearance 
and propagation of the norm for solution of~\eqref{eq:base}, and thus which can be connected to the 
quantitative nonlinear results in~\cite{MoVi}. It leads to study 
the linearized operator $\ron{L}$ in the space $L^1(\exp ( a|v|^s))$, 
which has no hilbertian self-adjointness structure. 
\medskip

The main result of this work is:
\smallskip

  \begin{thm} \label{theo:cvg}
  Let $B$ be a collision kernel satisfying assumption~\eqref{eq:HS}. 
  Let $\lambda \in (0,\nu_0)$ be the spectral gap of the linearized operator $\ron{L}$ in $L^2(M^{-1})$.   
  Let $f_0$ be a nonnegative initial datum in $L^1(1+|v|^2)$. 
  Then the solution $f=f(t,v)$ to the spatially 
  homogeneous Boltzmann equation~\eqref{eq:base} with initial datum $f_0$ 
  satisfies: for any $0<\mu \le \lambda$,  
  there is a constant $C$, which depends explicitly on $B$, the mass and energy  
  of $f_0$, on $\mu$ and on a lower bound on $\nu_0 -\mu$, such that 
    \begin{equation*}\label{eq:cvg}
    \|f(t,\cdot)-M\|_{L^1(\R^N)} \le C \, e^{- \mu t}.
    \end{equation*}
  \end{thm}
\medskip

\begin{rem} 

%1. Note that the optimal rate $\mu=\lambda$ is allowed in the theorem, 
%which can be related to the fact that the eigenspace of 
%$\ron{L}$ associated with the first non-zero eigenvalue $-\lambda$ is not degenerate. 
%It seems to be the first time this optimal rate is reached, 
%since both the quantitative study in \cite{CaGaTo} 
%for Maxwell molecules and the non-constructive results of \cite{Ar88} for hard spheres 
%only prove a convergence like $O(e^{-\mu t})$ for any $\mu <\lambda$, where 
%$\lambda$ is the corresponding spectral gap. 
%\smallskip

Let us recall that from Section~\ref{sec:BM} (work~\cite{BaMo}),
one deduces (in dimension $3$ for $C_B =1$) 
that $\lambda$ is controlled from below by some $c \approx 0.04$. 
Hence by choosing $\mu =c$ we get a completely explicit rate of exponential 
convergence, since 
  \[ \nu_0 = |\ens{S}^2| \, \int_{\R^3} M(v) |v| \, dv \] 
is explicit. 
\end{rem}
\medskip

We also state the functional analysis result used in the proof 
of Theorem~\ref{theo:cvg} and which has interest in itself. 
We consider $\bar{\ron{L}}$, that is the unbounded linearized collision operator on the space $L^1(\exp ( a|v|^s))$ with 
$a>0$ and $0<s<\gamma/2$ (with domain $\mbox{Dom}(\bar{\ron{L}}) = L^1 (\nu \exp ( a|v|^s))$), 
and $\ron{L}$, that is the unbounded self-adjoint linearized collision operator on $L^2(M^{-1})$ 
(with domain $L^2(\nu^2 M^{-1})$).
These operators are closed and we have
\medskip

 \begin{thm}\label{theo:spectre}
 Let $B$ be a collision kernel satisfying  
 assumption~\eqref{eq:HS}. Then the spectrum $\Sigma(\bar{\ron{L}})$ of $\bar{\ron{L}}$ 
 is equal to the spectrum $\Sigma(\ron{L})$ of $\ron{L}$. 
 Moreover $\bar{\ron{L}}$ and $\ron{L}$ have the same 
 eigenvectors associated to discrete eigenvalues. 
 \end{thm}
 \medskip

\begin{rems} 

1. This theorem essentially means that enlarging the functional space from 
$g \in L^2(M^{-1})$ to $g \in L^1 (\exp ( a|v|^s))$ 
does not yield new eigenvectors for the linearized collision operator.
\smallskip

2. It implies in particular that $\bar{\ron{L}}$ only has  
non-degenerate eigenspaces associated with its discrete 
eigenvalues, since this is true for the self-adjoint operator $\ron{L}$. 
This is related to the fact that the optimal  
convergence rate $\mu=\lambda$ is reached. 
Note that it seems to be the first time that such a result is obtained, 
since both the quantitative study in \cite{CaGaTo} 
for Maxwell molecules and the non-constructive results of \cite{Ar88} for hard spheres 
only prove a convergence like $O(e^{-\mu t})$ for any $\mu <\lambda$, where 
$\lambda$ is the corresponding spectral gap. 
%is exactly $C \, e^{-\lambda t}$ and 
%$not $C \, t^k \, e^{-\lambda t}$ for some $k >0$. 
It also yields a simple form of the 
first term in the asymptotic development (see \cite{GM:04}).
\end{rems}
\medskip

%3. This result means that, for hard spheres, 
%the linear part of the collision operator $f \to Q(M,f) + Q(f,M)$ 
%``has a spectral gap'' in $L^1(\exp ( a|v|^s))$, in the sense 
%that it satisfies exponential decay estimates 
%on its evolution semigroup in this space, with the 
%rate given by the spectral gap of $\ron{L}$. 
%We use this linear feature of the collision process to compensate for 
%the fact the functional inequality ({\em Cercignani's conjecture})
%  \begin{equation}\label{eq:CerCon}
%  {\cal D}(f) \ge K \, \big[ H(f) - H(M) \big], \quad K>0, 
%  \end{equation}
%is not true for $f \in L^1(1+|v|^2)$. It also supports 
%the fact that \eqref{eq:CerCon} could be true for 
%%that there is non ``nonlinear spectral gap'' in $L^1 _2$,  
%%i.e. Cercignani's conjecture on entropy production is not true in this space. 
%%It also supports the fact that Cercignani's conjecture could be true for 
%solutions $f$ of~\eqref{eq:base} satisfying some exponential decay at infinity 
%(as was questioned in~\cite[Chapter~3, Section~4.2]{Vi:hb}), in the sense $f \in L^1(\exp ( a|v|^s))$. 
%\medskip

Let us explain shortly the method of the proof. 
The main idea is to establish 
quantitative estimates of exponential decay on the evolution 
semigroup of $\bar{\ron{L}}$. They are used to estimate the rate of convergence 
when the solution is close to equilibrium (where the linear part of the collision 
operator is dominant), whereas the existing nonlinear entropy method, combined  
with some {\em a priori} estimates in $L^1(\exp ( a|v|^s))$, is used to estimate 
the time necessary to enter a given neighborhood of equilibrium, 
\smallskip

Concerning the first step, that is to show the exponential decay 
on the semigroup of $\bar{\ron{L}}$, it is done by proving that $\bar{\ron{L}}$ is 
sectorial, with explicit control on the resolvent, and with the same spectrum as $\ron{L}$ (see the introduction 
for the definition of a sectorial operator). 
This relies essentially on Hilbert's decomposition of $\bar{\ron{L}}$ between a local part 
(the multiplication by the collision frequency) and a non-local 
part, and then on the proof that the non-local part is well approximated 
by some truncated operator mapping functions of $L^1(\exp ( a|v|^s))$ 
into $L^2(M^{-1})$. From this one can deduce that the resolvent 
of $\bar{\ron{L}}$ at $\xi \in \C$ exists as soon as the one 
of $\ron{L}$ exists, with an explicit relation between their norms. 
This shows that $\Sigma(\bar{\ron{L}}) \subset \Sigma(\ron{L})$, and yields
the desired sectoriality estimate on $\bar{\ron{L}}$. More accurate  
results on the essential spectrum of $\bar{\ron{L}}$ (in order 
to obtain Theorem~\ref{theo:spectre}) can be obtained 
using perturbative arguments as Grad did for $\ron{L}$. 
\smallskip

Concerning the second step, that is the connection to the nonlinear theory, 
it is mainly based on three ingredients: a Gronwall argument to deal 
with the remaining quadratic term, some sharp versions of Povzner 
inequalities (in the spirit of \cite{Bo:mts,BoGaPa}) to show the appearance 
and propagation of the $L^1(\exp ( a|v|^s))$ norm, and the quantitative result 
of convergence to equilibrium of \cite{MoVi} to obtain an explicit time 
for the solution to the nonlinear Boltzmann equation to enter into some 
given neighborhood of equilibrium. 
\medskip

%{\bf Acknowledgment.}
%\medskip

\end{document}